\documentclass[a4paper,12pt]{article}
\usepackage{amsmath,amssymb} \usepackage{euler,eucal}
\usepackage[utf8]{inputenc}
\usepackage[T1]{fontenc}
\usepackage[english]{babel}
\usepackage[unicode=true,plainpages=false]{hyperref}
\usepackage{graphicx}
\hypersetup{colorlinks=false,pdfborder={0 0 0.1},linkcolor=magenta,anchorcolor=magenta,urlcolor=blue,citecolor=blue,
          pdftitle={A tensor decomposition algorithm for large ODEs with conservation laws},
          pdfauthor={S. V. Dolgov},
          pdfkeywords={high--dimensional problems, tensor train format, DMRG, alternating iteration, differential equations, conservation laws}
          }

\usepackage[left=30mm,right=30mm,top=20mm,bottom=30mm,a4paper]{geometry}
\usepackage{algpseudocode}
\usepackage[ruled]{algorithm}
\usepackage{amsthm}
\theoremstyle{definition}

\newtheorem{remark}{Remark}

\usepackage{tikz} \usetikzlibrary{positioning,shapes}
\usepackage{pgfplots} \usetikzlibrary{plotmarks}
\usetikzlibrary{calc,backgrounds}

\input{dats/pgfart.sty}

\def\O{\mathcal{O}}
\def\eps{\varepsilon}
\def\C{\mathbb{C}}

\def\x{\mathbf{x}}

\def\z{\mathbf{z}}
\def\B{\mathbf{B}}

\def\I{\mathcal{J}}
\def\le{\leqslant}
\def\leq{\leqslant}
\def\ge{\geqslant}

\date{September 30, 2014}

\title{A tensor decomposition algorithm for large ODEs with conservation laws%
\thanks{The author acknowledges funding from the EPSRC fellowship EP/M019004/1.%
         }}
\author{Sergey V. Dolgov\footnote{University of Bath, Claverton Down, BA2 7AY, Bath, United Kingdom ({\tt s.dolgov@bath.ac.uk}).} }

\begin{document}

\maketitle

\begin{abstract}
We propose an algorithm for solution of high-dimensional evolutionary equations (ODEs and discretized time-dependent PDEs) in the Tensor Train (TT) decomposition, assuming that the solution and the right-hand side of the ODE admit such a decomposition with a low storage.
A linear ODE, discretized via one-step or Chebyshev differentiation schemes, turns into a large linear system.
The tensor decomposition allows to solve this system for several time points simultaneously using an extension of the Alternating Least Squares algorithm.
This method computes the TT approximation of the solution directly, without ever solving the original large problem,
and encapsulates the Galerkin model reduction of the ODE.
This allows an efficient estimation of the time discretization error, and hence provides a way to adapt the time steps.
Besides, conservation laws can be preserved exactly in the reduced model by expanding the approximation subspace with the generating vectors of the linear invariants and correction of the euclidean norm.
In numerical experiments with the transport and the chemical master equations,
we demonstrate that the new method is faster than traditional time stepping and stochastic simulation algorithms, whereas the invariants are preserved up to the machine precision irrespectively of the TT approximation accuracy.

\end{abstract}

\section{Introduction}
Large-scale evolutionary equations for many-body systems arise ubiquitously in the numerical modeling.
The cases of particular interest and difficulty involve many configuration coordinates in the state space.
For instance, the time-dependent \emph{Schroedinger} equation describes the wavefunction, which depends on all positions of all quantum particles or states of spins.
Another important example is the joint probability density function, which obeys \emph{Fokker-Planck} or \emph{master} equations in continuous or discrete spaces, respectively.
The solution of a problem with $d$ configuration variables is a $d$-variate function.
When $d$ is much larger than $3$, a reasonably uniform discretization would require $\O(n^d)$ degrees of freedom.
Typical examples in quantum physics involve $d$ being of the order of hundreds, and the straightforward computation with $n^d$ unknowns is impossible.

To cope with such \emph{high-dimensional} problems, one has to employ \emph{data-sparse} techniques, i.e. describe the solution by much fewer unknowns than $n^d$.
Different approaches exist for this task.
Among the most successful ones we may identify Monte Carlo (and Quasi Monte Carlo) methods \cite{nieder-qmc-1978,graham-QMC-2011}, Sparse Grids \cite{smolyak-1963,griebel-sparsegrids-2004}, and tensor product representations.
In this paper, we adopt the latter framework.

\emph{Tensor product decompositions} rely on the idea of separation of variables: a $d$-variate array (or \emph{tensor}) can be defined or approximated by sums of products of univariate factors.
Extensive information can be found in recent reviews and books, e.g. \cite{kolda-review-2009,hackbusch-2012,bokh-surv-2015}.
A promising potential of the tensor product methods stems from the fact that a univariate factor is defined by only $n$ values.
If a tensor can be approximated up to the required accuracy with a moderate number of factors, the memory and complexity savings can be outstanding.

There exist different tensor product \emph{formats}, i.e. rules that map univariate factors to the initial array.
In case of two dimensions, one ends up with the low-rank dyadic factorization of a matrix.
This straightforward sum of direct products of vectors in higher dimensions is called the CP format \cite{hitchcock-rank-1927}.
However, the CP approximation problem may be ill-posed \cite{desilva-2008}.
This problem is circumvented in recurrent two-dimensional factorizations, where one can enforce a certain stable form of the representation.
In this paper, we focus on the simplest example,
the so-called \emph{Tensor Train} (TT) decomposition \cite{osel-tt-2011}.
It was rediscovered several times, and the most important analogs in quantum physics are \emph{Matrix Product States} (MPS) \cite{fannes-mps-1992} and
\emph{Density Matrix Renormalization Group} (DMRG) \cite{white-dmrg-1993}.
This format possesses all the power of recurrent factorizations, but algorithms are easier to describe.
For higher flexibility in particular problems, one can use more general tree-based constructions, such as \emph{HT} \cite{gras-hsvd-2010} or \emph{Extended TT/QTT-Tucker} \cite{dk-qtt-tucker-2013} formats.

DMRG is not only the name of the representation, but also a variety of computational tools.
It was originally developed for finding ground states (lowest eigenpairs) of high-dimensional Hamiltonians of spin chains.
The main idea behind DMRG is the alternating optimization of a function (e.g. Rayleigh quotient) over tensor format factors.
It was noticed that this method may manifest a remarkably fast convergence,
and later extensions to the energy function followed \cite{jeckelmann-dmrgsolve-2002,holtz-ALS-DMRG-2012}.

Besides the stationary problems, the same framework was applied to the dynamical spin Schroedinger equation.
Two conceptually similar techniques, the \emph{time-evolving block decimation} (TEBD) \cite{vidal-tebd-2004} and the \emph{time-dependent DMRG} (tDMRG) \cite{white-tdmrg-2004} take into account the nearest-neighbor form of the Hamiltonian to split the operator exponent into two parts using the Trotter decompositions.
Each part can then be integrated exactly, followed by the separation of variables via the truncated singular value decomposition.
This methods performs well for short times, but in a long time integration the error may accumulate, and the storage of the tensor product decomposition grows dramatically \cite{schollwock-2011}.

To avoid this problem, one can use the so-called \emph{Dirac-Frenkel} principle \cite{lubich-koch-dynten-2010,lob-ttksl-2015}.
This scheme projects the dynamical equations onto the tangent space of the tensor product manifold.
The storage of the format is now fixed, but the approximation errors can be difficult to control.

As an alternative approach, we consider time as an extra variable and solve a global system for many time steps simultaneously \cite{PeterSchwab-parabolic-2004,DKhOs-parabolic1-2012,kazeev-hp-dg-time-2012}.
Having several time steps allows to estimate the time discretization error and adapt the time grid accordingly.
However, the global state-time system is non-symmetric and requires a reliable solution algorithm in the tensor format.
We use an extension of DMRG, the so-called \emph{Alternating Minimal Energy} (AMEn) method \cite{ds-amen-2014}. It augments the tensor format of the solution by the tensor format of the global residual.
This improves the convergence and allows to adapt the tensor format storage up to a desired accuracy tolerance.

The residual is not the only quantity we can enrich the solution with.
The approximation error of the tensor decomposition is distributed evenly in all components of the solution.
However, it might be beneficial to compute some parts of the solution with a higher accuracy.
For example, the exact ODE may possess certain \emph{conservation laws} (e.g. phase \cite{sdwk-nmr-2014} or normalization), which are worth to be preserved in a numerical scheme.
We show that the basis vectors of the co-kernel of the ODE matrix can be inserted into the TT representation of the solution in addition to the residual.
This allows to preserve the corresponding linear invariants.
The second norm of the solution can then be corrected by rescaling.

The paper is structured as follows.
In the next section we formulate the ODE problem, 
investigate its properties related to the first- and the second-order invariants,
show the Galerkin model reduction concept and how the invariants can be preserved in the reduced system,
and suggest an adaptive linear discretization in time.
Section \ref{sec:tensor} starts with a brief introduction to tensor product formats and methods
and presents the new tAMEn algorithm (the name is motivated by tDMRG).
Section 4 demonstrates supporting numerical examples, followed by the conclusion in Section 5.

\section{Ordinary differential equations}
Our problem of interest is a linear system of ODEs,
\begin{equation}
\frac{d x}{dt}  = A(t)x, \qquad x(0)  = x_0,
\label{eq:ode}
\end{equation}
solved on $t \in [0,T]$, where $A(t) \in \C^{N\times N}$ is a stable matrix.
Throughout the paper, $x$ and other quantities denoted by small letters are $N \times 1$ vectors, such that the inner products can be consistently written as $c^* x \in \C^{1 \times 1}$.
The formulation \eqref{eq:ode} can be extended to ODEs with forcing, $dx/dt=Ax+f(t)$, or weakly nonlinear systems,
where $A(t)=A(t,\check{x}(t))$ depends on the solution from the previous Picard iteration.

\subsection{Conservation laws and Galerkin reduction}
Our goal will be to seek an ODE solution in a compressed data-sparse form.
A particular question of interest is the following: if the system preserves some quantities in time, is it possible to maintain this property in data-sparse algorithms, which are based on the Galerkin projection approach?

The simplest conservation laws are defined by linear functions of the solution and its euclidean norm.
Given some \emph{detecting} vector $c \neq 0$, the linear function can be written as $c^* x$.
It corresponds, for example, to the probability normalization in the Fokker-Planck equation: $x$ represents the discretized probability density function, and $\sum_{i=1}^{N} x(i) = c^*x=1$, with $c$ being a vector of all ones.
For a time-invariant system $dx/dt=Ax$, a sufficient condition for conservation of $c^*x$ is the nullspace equation $A^* c = 0$.

Among the second-order invariants, we consider the euclidean (Frobenius) norm of the solution, $\|x\| = \sqrt{x^*x}$.
The conservation law $\|x(t)\|=\|x_0\|$ is a well-known property of the Schroedinger equation $dx/dt = \mathrm{i} H x$, where $\mathrm{i}$ is the imaginary unity, and $H=H^\top \in \mathbb{R}^{N\times N}$.
A sufficient condition is the skew-symmetry of the matrix, $A=-A^*$.

An abstract Galerkin reduction can be written as follows.
Given an orthogonal set of columns $X \in \C^{N \times r}$, $X^*X=I$, we replace the large system \eqref{eq:ode} by a reduced ODE\footnote{for simplicity, we consider the time-invariant ODE in this section.},
\begin{equation}
\frac{d v}{dt}  = \left(X^*AX\right) v, \qquad v(0)  = v_0= X^* x_0.
\label{eq:ode_red}
\end{equation}
Numerical treatment of this equation is cheap if the basis size is small, $r \ll N$.
The approximate solution of the initial problem \eqref{eq:ode} writes as $\tilde x(t) = X v(t)\approx x(t)$.
Many approaches exist to determine the basis sets $X$, see e.g. \cite{antoulas-modredsurvey-2001,benner-modredsurvey-2014}.
The Krylov method for the computation of the matrix exponential \cite{moler-expm-2003} belongs to this class as well.
Another celebrated technique is the Proper Orthogonal Decomposition (POD) \cite{lumley-turbulent-1967,sirovich-turbulence-1987,kerschen-pod-2005,nouy-modred-2010},
which extracts principal components from a set of \emph{snapshots} $\{x(t_j)\}_{j=1}^{\I}$ using the singular value decomposition.

The accuracy $\|x-\tilde x\|$ of the reduced model depends on the approximation capacity of the basis set.
In this paper, we employ a tensor product algorithm, which is similar to POD but computes both the basis and the reduced solution adaptively without solving the large original problem.
Most importantly, it belongs to the Galerkin projection framework \eqref{eq:ode_red}.
Here we show how to preserve first and second order invariants with an arbitrary Galerkin basis.

Suppose we are given vectors $C = \begin{bmatrix}c_1 & \cdots & c_M\end{bmatrix}$ such that $A^* C=0$.
Let us \emph{expand} the basis by concatenating $C$ and $X$ and orthogonalizing the result,
\begin{equation}
\begin{bmatrix}C & X\end{bmatrix} = \hat X R,  \quad \hat X^* \hat X = I \quad \mbox{(QR decomposition)}.
 \label{eq:X_enrich_C}
\end{equation}
Since the first $M$ columns of $\hat X$ belong to the kernel of $A^*$, the reduced matrix writes
$$
\hat X^* A \hat X =  \begin{bmatrix}\mathcal{C}^* A \mathcal{C} & \mathcal{C}^* A \mathcal{X} \\ \mathcal{X}^* A \mathcal{C} & \mathcal{X}^* A \mathcal{X} \end{bmatrix} = \begin{bmatrix}0 & 0 \\ \mathcal{X}^* A \mathcal{C} & \mathcal{X}^* A \mathcal{X}\end{bmatrix}, \quad \mbox{where} \quad \hat X = \begin{bmatrix}\mathcal{C} & \mathcal{X}\end{bmatrix}.
$$

In order to derive the reduced solution $v(t) = \exp\left(t \hat X^* A \hat X\right) v_0$ in the expanded basis,
consider one recursion step for the exponential series. For any $k=1,2,\ldots,$
$$
\begin{bmatrix}0 & 0 \\ \left(\mathcal{X}^* A \mathcal{X}\right)^{k-1} \mathcal{X}^* A \mathcal{C} & \left(\mathcal{X}^* A \mathcal{X}\right)^k\end{bmatrix} \begin{bmatrix}0 & 0 \\ \mathcal{X}^* A \mathcal{C} & \mathcal{X}^* A \mathcal{X}\end{bmatrix} = \begin{bmatrix}0 & 0 \\ \left(\mathcal{X}^* A \mathcal{X}\right)^{k} \mathcal{X}^* A \mathcal{C} & \left(\mathcal{X}^* A \mathcal{X}\right)^{k+1}\end{bmatrix},
$$
and hence we obtain
\begin{equation}
\exp\left(t \hat X^* A \hat X\right) = I + \sum_{k=1}^{\infty} \frac{\left(t \hat X^* A \hat X\right)^k}{k!} =  \begin{bmatrix}I & 0 \\  \sum\limits_{k=1}^{\infty} \frac{t \left(t \mathcal{X}^* A \mathcal{X}\right)^{k-1}}{k!} \mathcal{X}^* A \mathcal{C} & \exp\left(t \mathcal{X}^* A \mathcal{X}\right)\end{bmatrix}.
\end{equation}
Since the first row contains only the identity, the linear invariants $\mathcal{C}^* x_0$ are explicitly preserved in the solution,
$
v(t) = \begin{bmatrix}\mathcal{C}^* x_0 \\ w(t)\end{bmatrix}.
$

The skew-symmetry, yielding conservation of the second norm, is even easier to consider, since for any Galerkin projection,
$\left(X^* A X \right)^* = X^* A^* X = -X^* A X$, and hence $\|v(t)\| = \|X^* x_0\|$.
Moreover, $\|\tilde x(t)\| = \|v(t)\|=\|X^* x_0\|$ if $X$ is orthogonal.
Thus, it is enough to guarantee $\|X^* x_0\|=\|x_0\|$.
A simple way to do this is to rescale the projected initial state.
However, this requires a certain care if we need to preserve both the second norm and the linear invariants.
Given $v_0 = \begin{bmatrix}\mathcal{C}^* x_0 \\ \mathcal{X}^* x_0\end{bmatrix}$, 
we can rescale only the bottom part.
This means finding $\theta>0$ such that
\begin{equation}
\|\hat v_0\|^2 =  \|\mathcal{C}^* x_0\|^2 + \theta^2 \|\mathcal{X}^* x_0\| = \|x_0\|^2, \quad \mbox{hence} \quad \theta = \frac{\sqrt{\|x_0\|^2 - \|\mathcal{C}^* x_0\|^2}}{\|\mathcal{X}^* x_0\|},
\label{eq:2norm-rescale}
\end{equation}
and the rescaled initial state reads
$$
\hat v_0 = \begin{bmatrix}\mathcal{C}^* x_0 \\ \theta \mathcal{X}^* x_0\end{bmatrix}.
$$
Due to orthogonality of $\mathcal{C}$ and $\mathcal{X}$, it holds that $\|\mathcal{C}^* x_0\|\le \|\hat X^* x_0\|\le \|x_0\|$,
and hence $\theta$ is well-defined when $x_0 \notin \mathrm{span}(\mathcal{C})$.
Otherwise, $\|\mathcal{X}^* x_0\|=0$, and $\theta$ can be arbitrary.

\subsection{Linear discretization in time}\label{sec:spectral}
Assuming the solution $x(t)$ to be continuous, we can introducing a time discretization grid $\mathbf{t}=\{t_j\}_{j=1}^{\I} \in [0,T]$ and
collocate the solution on this grid, $\{x_j\}=\{x(t_j)\}$.
An approximate solution at any time can be computed by polynomial interpolation,
\begin{equation}
x(t) \approx \sum_{j=1}^{\I} x_j p_j(t),
\label{eq:tinterp}
\end{equation}
where $p_j(t)$ are polynomials, centered at $t_j$, such as the global Lagrange polynomials or local splines.
Since both ODE \eqref{eq:ode} and the interpolation \eqref{eq:tinterp} is linear in $x$, the discrete system is linear as well, and can be generally written as
\begin{equation}
Bx = f, \qquad B=I_N \otimes S - \left(I_N \otimes P \right) A(\mathbf{t}), \quad f =  x_0\otimes(Se),
\label{eq:xt-system}
\end{equation}
where $A(\mathbf{t})$ is a block-diagonal matrix constructed from the ODE matrices at the grid points,
\begin{equation}
A(\mathbf{t}) = \begin{bmatrix}A(t_1) \\ & \ddots \\ & & A(t_{\I})\end{bmatrix}, \qquad x = \begin{bmatrix}x(t_1) \\ \vdots \\ x(t_{\I})\end{bmatrix}
\label{eq:Adiag}
\end{equation}
is the vector of all snapshots stacked together,
$S\in\mathbb{R}^{\I \times \I}$ is the stiffness matrix corresponding to the time derivative,
$P\in\mathbb{R}^{\I \times \I}$ is the mass matrix,
$e = (1,\ldots,1)^\top\in\mathbb{R}^{\I}$ is a vector of all ones,
and $\otimes$ is the Kronecker product.
For a time-invariant ODE \eqref{eq:Adiag} simplifies to $A(\mathbf{t}) = A\otimes I_{\I}$.
For example, Euler and Crank-Nicolson schemes belong to this class with $S = \mathrm{tridiag}(-1,1,0)$,
and $P=\frac{T}{\I} I$ for the implicit Euler scheme on a grid $t_j = T j/\I$, and
$$
P = \frac{T}{2(\I-1)} \begin{bmatrix}
                       0 \\
                       1 & 1 \\
                       & \ddots & \ddots \\
                       & & 1 & 1
                      \end{bmatrix}
$$
for the Crank-Nicolson scheme on a grid $t_j = T(j-1)/(\I-1)$.
With these schemes we can take linear splines at $t_j \le t \le t_{j+1}$ in the interpolation \eqref{eq:tinterp}.

Alternatively, we can choose the Chebyshev grid with the nodes $t_j = \frac{T}{2}(1-\cos(\pi j/\I))$,
and use the Chebyshev spectral differentiation matrix \cite[Chapter 6]{trefethen-spectral-2000} $S = \{dp_j(t_i)/dt\}_{i,j=1}^{\I}$ , where $p_j$ is the Lagrange polynomial centered at $t_j$, and $P=I$.
An advantage of the spectral discretization is a rapid convergence (exponential in $\I$, see \cite{Tadmor-exp_acc_diff-1986} and \cite[Theorem 6]{trefethen-spectral-2000}) when the solution is analytic on the Bernstein ellipse extension of $[0,T]$.
On the other hand, lower order schemes lead to sparse matrices and lower condition numbers in \eqref{eq:xt-system}.

The Galerkin reduction \eqref{eq:ode_red} can be combined with \eqref{eq:xt-system} straightforwardly.
Given an orthogonal basis matrix $X$, we assemble and solve the $r\I \times r\I$ system
\begin{equation}
\left(I_r \otimes S - \left(I_r \otimes P\right)(X \otimes I_\I)^* A(\mathbf{t}) (X \otimes I_\I) \right) v =   v_0 \otimes (Se),
\label{eq:xt-system-red}
\end{equation}
where $v_0 = X^*x_0$.
Both linear and quadratic invariants can be preserved as shown in \eqref{eq:X_enrich_C} and \eqref{eq:2norm-rescale}, respectively.

\begin{remark}
Low-order schemes are often preferred to the spectral discretization because of the particular sparsity of the stiffness and mass matrices, e.g. \emph{bidiagonality}, which allows to solve \eqref{eq:xt-system} step by step.
However, in this paper we solve \eqref{eq:xt-system} indirectly via iterative tensor product algorithms (see Sec. \ref{sec:als}),
which require a single system of equations, defining the entire solution.
On the other hand, tensor decompositions allow more freedom in the choice of $S$ and $P$ due to the reduced cost; in fact, solving the global system \eqref{eq:xt-system} can be faster and more accurate than the step by step integration \cite{DKhOs-parabolic1-2012}, since it allows to take more accurate time discretization.
\end{remark}

\begin{remark}
If the ODE solution lacks smoothness, more sophisticated Discontinuous Galerkin techniques may be required \cite{schotzau-hp-dg-time-1999,kazeev-hp-dg-time-2012}.
Otherwise, the collocation leads to easier pointwise construction of the matrix \eqref{eq:Adiag}, compared to the computation of the Galerkin coefficients.
\end{remark}

An analog of the Runge's rule \cite{hall-ode-1976} can be used for estimating the discretization error.
Consider two grids with $\I$ and $2\I$ points, $\{t_j\}_{j=1}^{\I}$ and $\{t_i^*\}_{i=1}^{2\I}$.
Given an approximation $y(t) \approx dx/dt$ on the coarse grid $\{t_j\}$ (in our case $y(t) = A(t)x(t)$),
we can take the difference on the fine grid $|dx/dt(t_i^*) - y(t_i^*)|$ as our error estimate.
For evaluating the quantities on $\{t_i^*\}$ we construct
the fine-grid differentiation matrix $\hat S \in \mathbb{R}^{2\I \times 2\I}$
and the interpolation matrix $\hat P\in\mathbb{R}^{2\I \times \I}$,
which maps from $\{t_j\}$ to $\{t_i^*\}$.
Then the estimate can be computed from the snapshots as follows,
\begin{equation}
\mathcal{E}_{\I,T} = \left\|\left[I_N \otimes (\hat S \hat P) - (I_N \otimes \hat P)A(\mathbf{t})\right]x -  x_0 \otimes (\hat S \hat e) \right\|,
\label{eq:time_res}
\end{equation}
where $\hat e$ is a vector of all ones of size $2\I$.
For the Chebyshev discretization, for example, $\hat P_{i,j} = p_j(t_i^*)$.

\section{Tensor product representations and methods}\label{sec:tensor}
\subsection{Vectors and tensors}
The unknowns in the whole discrete solution can be enumerated by at least two independent indices, corresponding to the state space and time.
Assuming that $i=1,\ldots,N$ enumerates the state components of the solution, $x_i(t)$, and that the time points are enumerated by an index $j=1,\ldots,\I$, we can consider the solution as a matrix $X = \left[x_i(t_j)\right]$.
Moreover, we will assume (and exploit) that the state space can be further factorised into $d$ independent indices $i_1,\ldots,i_d$, running from $1$ to $n_1,\ldots,n_d$, respectively.
An equivalence between \emph{digits} $i_1,\ldots,i_d$ and the original index $i$ holds due to the standard positional expression,
\begin{equation}
i = (i_1-1)n_2\cdots n_d + (i_2-1) n_3 \cdots n_d + \cdots + i_d.
\label{eq:multiind-i}
\end{equation}
However, the solution can now be also seen as a \emph{tensor},
$\x=\left[\x(i_1,\ldots,i_d,j)\right] \in\C^{n_1\times\cdots\times n_d \times \I}$.
The multi-index expansion can arise for example from a discretization of PDEs:
if a PDE $\frac{\partial x}{\partial t}(q_1,\ldots,q_d,t) = Ax(q_1,\ldots,q_d,t)$
is discretized in $q_1,\ldots,q_d$ by collocation on a Cartesian product of independent univariate grids $\{q_k(i_k)\}$, $k=1,\ldots,d$,
the nodal values of $x$ can be collected into a tensor $\x$, as described above.

To write the global state-time system \eqref{eq:xt-system} consistently, we need to reshape the whole tensor $\x$ into a vector $x$ of size $(n_1 \cdots n_d) \I \times 1$.
We can extend \eqref{eq:multiind-i} to any set of indices, introducing a general \emph{multi-index}
\begin{equation}
\overline{i_p \ldots i_q} = (i_p-1)n_{p+1} \cdots n_q + \cdots + i_q, \quad q\ge p.
\label{eq:multiind}
\end{equation}
Now we can address the solution by either of the equivalent forms $x(\overline{i_1\ldots i_d,j})$, $X(\overline{i_1\ldots i_d},j)$ or $\x(i_1,\ldots,i_d,j)$.

\subsection{Tensor Train decomposition}\label{sec:tt}
The Tensor Train (TT) \cite{osel-tt-2011}, or Matrix Product States (MPS) \cite{fannes-mps-1992} decomposition for the tensor $\x$ (resp. vector $x$) is defined as follows,
\begin{equation}\label{eq:tt}
  x(\overline{i_1\ldots i_d,j})  =
  \sum_{\alpha_1=1}^{r_1} \cdots \sum_{\alpha_{d}=1}^{r_{d}} \x^{(1)}_{\alpha_1}(i_1) \x^{(2)}_{\alpha_1,\alpha_2}(i_2) \cdots \x^{(d)}_{\alpha_{d-1},\alpha_d}(i_d) \x^{(d+1)}_{\alpha_d}(j).
\end{equation}
The summation indices $\alpha_k=1,\ldots,r_k$, $k=1,\ldots,d$, are called the \emph{rank} indices,
and their ranges $r_k$ are the \emph{tensor train} ranks (TT ranks).
The right-hand side consists of the TT \emph{blocks} $\x^{(k)}\in\C^{r_{k-1}\times n_k\times r_k}$.
Introducing uniform bounds $r_k \le r$, $n_k \le n$, we can estimate the storage complexity of the TT decomposition as $\O(dnr^2)$.
If the rank bound $r$ is small, this is much lower than $N \I = \O(n^d \I)$ in the straightforward storage of $x$.

The matrix $B$ from \eqref{eq:xt-system} can be seen as a $(2d+2)$-dimensional tensor
and decomposed in a slightly different \emph{matrix} TT decomposition,
\begin{equation}
 B(\overline{i_1\ldots i_d,j},\overline{i_1\ldots i_d,j}')
        = \sum_{\gamma_1=1}^{R_1} \cdots \sum_{\gamma_{d}=1}^{R_{d}} \B^{(1)}_{\gamma_1}(i_1,i_1') \cdots \B^{(d)}_{\gamma_{d-1},\gamma_d}(i_d,i_d') \B^{(d+1)}_{\gamma_d}(j,j').
\label{eq:ttm}
\end{equation}
The matrix TT decomposition is introduced for consistency with the Kronecker product when $R_1=\cdots=R_d=1$
and multiplication with a ``vector'' TT decomposition of $x$ \eqref{eq:tt}.

The multi-index notation \eqref{eq:multiind} allows to notice that the TT decomposition can be seen as a low-rank decomposition of a
matrix $X^{\{k\}} = \left[X(\overline{i_1\ldots i_k}, \overline{i_{k+1}\ldots j})\right]$ for any $k=1,\ldots,d$.
We can group the left, respectively, right subset of TT blocks into \emph{interface matrices}, or simply \emph{interfaces}
\begin{equation}\label{eq:iface}
\begin{split}
 X^{(\le k)}(\overline{i_1\ldots i_k}, \alpha_k) & = \sum_{\alpha_1=1}^{r_1} \cdots \sum_{\alpha_{k-1}=1}^{r_{k-1}} \x^{(1)}_{\alpha_1}(i_1)  \cdots \x^{(k)}_{\alpha_{k-1},\alpha_k}(i_k), \\
 X^{(>k)}(\alpha_k, \overline{i_{k+1}\ldots j}) & = \sum_{\alpha_{k+1}=1}^{r_{k+1}} \cdots \sum_{\alpha_{d}=1}^{r_{d}} \x^{(k+1)}_{\alpha_k,\alpha_{k+1}}(i_{k+1}) \cdots \x^{(d+1)}_{\alpha_{d}}(j).
\end{split}
\end{equation}
We can naturally extend this definition to $X^{(<k)} = X^{(\le k-1)}$ and $X^{(\ge k)} = X^{(>k-1)}$.
Then we can write $X^{\{k\}} = X^{(\le k)} X^{(>k)}$.
Moreover, the interface matrices allow to see the TT decomposition as a \emph{linear map} of each TT block $\x^{(k)}$.
Indeed, reshaping it into a vector $x^{(k)}(\overline{\alpha_{k-1}i_k\alpha_k}) = \x_{\alpha_{k-1},\alpha_k}(i_k)$,
we can write $x = X_{\neq k} x^{(k)}$, where
$X_{\neq k}$ is the \emph{frame} matrix
\begin{equation}\label{eq:frame}
 X_{\neq k}
  = X^{(<k)} \otimes I_{n_k} \otimes \left(X^{(>k)}\right)^\top.
\end{equation}

\subsection{Computing TT decompositions by alternating iteration}\label{sec:als}
Although a TT approximation can be computed for any tensor via a sequence of singular value decompositions (SVD) \cite{osel-tt-2011},
it is rarely efficient or even possible when the tensor is large.
The aim of the tensor product methodology is to avoid fully stored tensors at all stages of computations.
One of the most successful approaches traces back to the alternating least squares optimization over the tensor decomposition blocks.
It was then generalized to the Alternating Linear Scheme (ALS) \cite{holtz-ALS-DMRG-2012}.
A similar algorithm, called Density Matrix Renormalization Group (DMRG) \cite{white-dmrg-1993,jeckelmann-dmrgsolve-2002}, was proposed in quantum physic for calculation of ground states, i.e. lowest eigenvalues of high-dimensional Hamiltonians.

Let us consider the linear system $Bx=f$ as an overdetermined equation on a particular block $x^{(k)}$ in the TT decomposition \eqref{eq:tt};
the linearity established in the previous subsection makes this equation linear, $(BX_{\neq k}) x^{(k)} = f$.
This equation can be resolved in different ways (e.g. by least squares), but practically the cheapest option is to use the frame matrix,
\begin{equation}
\left(X_{\neq k}^* B X_{\neq k}\right) x^{(k)} = X_{\neq k}^* f.
\label{eq:localsys}
\end{equation}
This reduction can be justified by relation to the minimization of the energy function $x^* B x - 2 \mathrm{Re}~x^* f$ when the matrix $B$ is symmetric positive definite (SPD).
However, the projection formalism \eqref{eq:localsys} is more general and can be applied also if $B$ is not SPD, which is the case for \eqref{eq:xt-system}.
The \emph{alternating} iteration is realised by sweeping through different blocks, $k=1,\ldots,d+1,$ and backwards from $k=d+1$ to $k=1$ until convergence, solving \eqref{eq:localsys} in each step.

Three essential details make the alternating iteration actually useful:
\begin{itemize}
 \item efficient assembly of \eqref{eq:localsys};
 \item orthogonality of $X_{\neq k}$ and efficient solution of \eqref{eq:localsys};
 \item adaptation of TT ranks of $x$.
\end{itemize}
The frame matrix, composed from the interface matrices \eqref{eq:iface} via Kronecker products, can be seen as a special TT decomposition with the same number of blocks as in $x$.
In turn, the matrix $B$ and right-hand side $f$ are assumed to be available in the TT format as well, such as \eqref{eq:ttm}.
This allows to compute $X_{\neq k}^* B X_{\neq k}$ and $X_{\neq k}^* f$ efficiently, using only the TT operations.
Moreover, sequential iteration over $k=1,2,\ldots$ allows to reuse partial products of the interfaces of $x$, $B$ and $f$ and make the algorithm even more efficient, with the total asymptotic complexity linear in $d$ \cite{holtz-ALS-DMRG-2012,DoOs-dmrg-solve-2011}.

The TT representation is not unique; any partition of identity can be inserted between adjacent TT blocks, e.g. $X^{\{k\}} = \left(X^{(\le k)} R\right)\left(R^{-1} X^{(>k)}\right)$, without changing the whole tensor.
However, the matrix $R$ changes the interfaces, and we can choose it in order to empower the representation with desirable properties.
For example, we can make $X^{(<k)}$ and $X^{(>k)}$ orthogonal by performing QR decompositions of appropriately reshaped TT blocks.
By construction \eqref{eq:frame}, $X_{\neq k}$ is orthogonal, too.
The orthogonality of the projection \eqref{eq:localsys} makes the reduced problem well conditioned,
which can be solved iteratively (we employ the BiCGstab algorithm) using fast matrix-vector products due to the TT structure inherited from the original problem \cite{DoOs-dmrg-solve-2011}.

For high-dimensional problems it is difficult to guess $d$ rank parameters.
It becomes necessary to adapt them during the computations in such a way that the TT solution is within the desired distance from the exact solution.
If we possess a solution with a satisfactory accuracy but overly large TT ranks, it is easy to reduce them via SVD \cite{osel-tt-2011}.
It is more important therefore to develop a procedure for increasing the ranks.
The DMRG method addresses this problem by reducing the system to a two-dimensional block (merged from $\x^{(k)}$ and $\x^{(k+1)}$), which can be split via SVD up to a desired threshold.
However, this requires solving a larger problem on the merged block.
The Alternating Minimal Energy (AMEn) algorithm \cite{ds-amen-2014} solves one-dimensional problems in each step,
but augments the TT blocks of the solution by the TT blocks of an approximate \emph{global} residual $z \approx f-Bx$.
Since $f,B$ and $x$ are all represented in the TT format, the residual can be approximated efficiently by the \emph{second} ALS iteration,
applied to a simpler problem $I z = f-Bx$.
Given a TT decomposition
$$
z(\overline{i_1\ldots i_d,j}) = \sum_{\beta_1,\ldots,\beta_d=1}^{\rho_1,\ldots,\rho_d} \z^{(1)}_{\beta_1}(i_1) \cdots \z^{(d+1)}_{\beta_d}(j)
$$
from the previous iteration,
we define the interface matrices $Z^{(<k)}$ and $Z^{(>k)}$ similarly to \eqref{eq:iface},
and update the $k$-th TT block of the residual by projecting
\begin{equation}
z^{(k)} = \left(Z^{(<k)} \otimes I \otimes (Z^{(>k)})^\top\right)^*(f-Bx).
\label{eq:zk}
\end{equation}
Performing this process simultaneously with the computation of the solution blocks \eqref{eq:localsys},
we ensure that $Z^{(<k)}$ and $Z^{(>k)}$ are sufficiently good bases for the residuals in all steps.
In turn, projecting the residual onto the solution interface,
\begin{equation}
\zeta^{(k)} = \left(X^{(<k)} \otimes I \otimes (Z^{(>k)})^\top\right)^*(f-Bx),
\label{eq:sk}
\end{equation}
we can expand the solution TT block,
\begin{equation}
 \x^{(k)}(i_k) = \begin{bmatrix}\x^{(k)}(i_k) & \boldsymbol\zeta^{(k)}(i_k)\end{bmatrix}, \qquad \x^{(k+1)}(i_{k+1}) = \begin{bmatrix}\x^{(k+1)}(i_{k+1}) \\ 0\end{bmatrix}.
\label{eq:enrich}
\end{equation}
This allows to increase the solution TT ranks (by the ranks of $\boldsymbol\zeta^{(k)}$),
and also improves convergence in difficult cases, since the basis of the reduction \eqref{eq:localsys} contains now the residual of the original problem.

\subsection{tAMEn: extended time integrator}
The time-dependent version of the AMEn algorithm combined two enrichments of the solution: by the residual \eqref{eq:enrich} and by the co-kernel vectors \eqref{eq:X_enrich_C}.
Assume the latter to be given in a compatible TT format,
$$
c_{m}(\overline{i_1 \ldots i_d}) = \sum_{\beta_1,\ldots,\beta_{d-1}=1}^{\rho_1,\ldots,\rho_{d-1}} \mathbf{c}^{(1)}_{\beta_1}(i_1) \mathbf{c}^{(2)}_{\beta_1,\beta_2}(i_2) \cdots \mathbf{c}^{(d)}_{\beta_{d-1},m}(i_d),
$$
where $m=1,\ldots,M$ enumerates different vectors $c_m$.
In the course of the alternating iteration from $k=1$ to $k=d$, the combined enrichment is performed as follows,
\begin{equation}
\x^{(k)}(i_k) = \begin{bmatrix}\x^{(k)}(i_k) & \boldsymbol\zeta^{(k)}(i_k) & \mathcal{C}_k \mathbf{c}^{(k)}(i_k) \end{bmatrix}, \qquad \x^{(k+1)}(i_{k+1}) = \begin{bmatrix}\x^{(k+1)}(i_{k+1}) \\ 0 \\ 0\end{bmatrix},
\label{eq:enrich_both}
\end{equation}
where $\mathcal{C}_k = (X^{(<k)})^* C^{(<k)}$ is the projection onto the left interface of the solution.
We can see that $c_m \in \mathrm{span}(X^{(<k)} \otimes I_{n_k \cdots n_d})$ for all $k=1,\ldots,d$ and $m=1,\ldots,M$.
For $k=1$, for example, we can write
$$
c_m = \left(\x^{(1)} \otimes I_{n_2 \cdots n_d}\right) \tilde c_m, \quad \tilde c_m = \begin{bmatrix}0 \\ 0 \\ c^{(>1)}_m\end{bmatrix},
$$
where $c^{(>1)}_m$ is the $(n_2\cdots n_d \rho_1) \times 1$ vectorisation of the interface matrix $C^{(>1)}_m$.
By induction, this extends to $k>1$.
In order to maintain orthogonality of the interfaces, we perform the QR decomposition of $\x^{(k)}$ after the enrichment \eqref{eq:enrich_both}.

For $k=d+1$, we can notice that the frame matrix reduces to $X^{(\le d)} \otimes I$.
Therefore, the local problem \eqref{eq:localsys} for $x^{(d+1)}$ is nothing else than the reduced discretized ODE \eqref{eq:xt-system-red}
with the interface being the Galerkin basis, $X = X^{(\le d)}$, and the last TT block being the unknown, $v=x^{(d+1)}$.
The enrichment \eqref{eq:enrich_both} ensures that this basis contains also the co-kernel matrix $C$.
Moreover, the second norm of the right hand side (reduced initial state) can be corrected according to \eqref{eq:2norm-rescale}.
If we stop the alternating iteration at this step, the error in the linear invariants and the second norm depends only on the accuracy of the solution of \eqref{eq:xt-system-red} and the time discretization, but not on the accuracy of the TT decomposition.
If the last TT rank $r_d$ is reasonably small, we can take sufficiently large $\I$ and solve \eqref{eq:xt-system-red} directly, which yields the machine precision accuracy in the conservation laws.

\begin{algorithm}[htb]
\caption{tAMEn algorithm}
\label{alg:tamen}%
\begin{algorithmic}[1]
 \Require Initial state $x_0$, matrix $A(t)$ and right hand side $f(t)$ in the TT format, final time $T$, accuracy threshold $\eps$, discretization points $\mathbf{t}\in [0,1]$ and matrices $S,P,\hat S$ and $\hat P$.
 \Ensure Time splitting points $T_0=0 < T_1 < \cdots < T_{L}=T$, solutions $x_{\ell}(t)$ in the TT format.
 \State Set $t=0$, $\ell=1$, $L=1$, $h_L=T$.
 \While{$t<T$}
   \State Rescale $\mathbf{t},S,P,\hat S$ and $\hat P$ from $[0,1]$ to $[T_{\ell-1}, T_{\ell-1}+h_\ell]$.
   \State Form $B = I \otimes S - (I \otimes P)\mathrm{diag}(A(\mathbf{t}))$ and $f = x_{\ell-1} \otimes (Se)$.
   \State Set $x=x_{\ell-1} \otimes e$.
   \For{$\mbox{iter}=1,2,\ldots,$}
     \State Set $x_{prev} = x$.
     \For{$k=d+1,d,\ldots,2$}
       \State Orthogonalize $X^{(>k)}$ and $Z^{(>k)}$, see \cite[Section 3]{osel-tt-2011}.
     \EndFor
     \For{$k=1,2,\ldots,d$} \Comment{Solve}
       \State Solve $(X_{\neq k}^* B X_{\neq k})x^{(k)} = X_{\neq k} f$, as defined in \eqref{eq:xt-system} and \eqref{eq:frame}.
       \State Compute truncated SVD of $\x^{(k)}$ up to $\eps$.
       \State Compute residual blocks as shown in \eqref{eq:zk} and \eqref{eq:sk}.
       \State Enrich $\mathbf{x}^{(k)}$ and $\mathbf{x}^{(k+1)}$ as shown in \eqref{eq:enrich_both}.
       \State Orthogonalize $X^{(<k+1)}$ and $Z^{(<k+1)}$, see \cite[Section 3]{osel-tt-2011}.
     \EndFor
     \State Correct the norm of $v_0 = (X^{(\le d)})^* x_{\ell-1}$ as shown in \eqref{eq:2norm-rescale}.
     \State Solve $(X_{\neq d+1}^* B X_{\neq d+1})x^{(d+1)} = \hat v_0 \otimes (S e)$.
     \State Compute the error estimate \eqref{eq:time_res_red} and $h_{\ell+1} =  h_{\ell} \cdot (\eps/\mathcal{E}_{\I,h_{\ell}})^{1/q}.$
     \If {$\mathcal{E}_{\I,h_{\ell}}\le\eps$}
       \If {$\|x-x_{prev}\|<\eps\|x\|$} \Comment{This step converged, accept it}
         \State Set $x_{\ell}=x$, $T_{\ell} = T_{\ell-1}+h_{\ell}$, $t=t+h_{\ell}$, $\ell=\ell+1$, and break.
       \EndIf
     \Else \Comment{Reject the step}
       \State Set $h_{\ell} = h_{\ell+1}$ and break.
     \EndIf
   \EndFor
 \EndWhile
\end{algorithmic}
\end{algorithm}


The time discretization error \eqref{eq:time_res} can be also estimated from the reduced system.
Instead of the full solution, we consider only the last TT block, and replace the state matrix by its projection\footnote{For non-autonomous ODEs the estimate can be extended accordingly.},
\begin{equation}
\mathcal{E}_{\I,T} = \left\|\left(I_{r_d} \otimes \hat S - \left[(X^{(\le d)})^* A  X^{(\le d)}\right] \otimes \hat P \right)x^{(d+1)} - \hat v_0 \otimes (\hat S e) \right\|.
\label{eq:time_res_red}
\end{equation}
This estimate can be used for refining the number of time points $\I$ or the length of the time interval.
Instead of solving \eqref{eq:xt-system} on the whole desired interval $[0,T]$,
we can split it into a sequence of subintervals $[0, T_1], \ldots, [T_{L-1}, T_L]$,
taking the solution at the last time point of the previous interval as the initial state in the next interval.
We determine an optimal splitting using the local error control with rejections \cite{byrne-ode-1975}.
We aim to maintain the error in the next time interval below a desired threshold,
$\mathcal{E}_{\I,h_{\ell+1}} \le \eps$,
so we adjust the next interval length as follows,
\begin{equation}
h_{\ell+1} = \mathcal{T}_{\ell+1} - \mathcal{T}_{\ell} = h_{\ell} \left(\dfrac{\eps}{\mathcal{E}_{\I,h_{\ell}}}\right)^{1/q}.
\label{eq:stepchange}
\end{equation}
The parameter $q$ reflects the order of convergence of the time scheme,
which is $1$ for the Euler method, $2$ for the Crank-Nicolson scheme, and $\I$ for the Chebyshev differentiation.
Moreover, if it appears that $\mathcal{E}_{\I,h_{\ell}} > \eps$, such solution is \emph{rejected}, the \emph{current} interval is shrunk according to \eqref{eq:stepchange}, and the solution is started again from $T_{\ell-1}$.
The entire procedure is written in Alg. \ref{alg:tamen}.
Assuming the index ranges from Sec. \ref{sec:tt} ($r_k \leq r$, $R_k \leq R$, $n_k \leq n$),
the computational complexity of Alg. \ref{alg:tamen}, inherited from AMEn \cite{ds-amen-2014}, reads
$$
\mathcal{O}(dn(Rr^3 + R^2 r^2)).
$$

\section{Numerical experiments}
We have implemented Algorithm \ref{alg:tamen} in Matlab, and carried out the computations on one core of the University of Bath \texttt{Balena} cluster, an Intel Xeon E5-2650 CPU at 2.6GHz.
The code is available from \href{http://github.com/dolgov/tamen}{\texttt{http://github.com/dolgov/tamen}}.

\subsection{Convection}
Our first example is the transport equation in the periodic domain $[-10,10]^2$ with the central difference discretization scheme,
\begin{equation}
\frac{dx}{dt} = \left(\nabla_n \otimes I_n + I_n \otimes \nabla_n\right) x, \quad \nabla_n = \frac{1}{2h} \begin{bmatrix}0 & 1 & \cdots & &  -1 \\ -1 & 0 & 1 \\ & \ddots & \ddots & \ddots \\ &  & -1 & 0 & 1 \\ 1 & & \cdots & -1 & 0 \end{bmatrix} \in \mathbb{R}^{n\times n},
\label{eq:conv}
\end{equation}
where $h = 20/n$ is the mesh step of the uniform grid $q_k(i_k) = -10+h (i_k-1)$, $i_k=1,\ldots,n$, $k=1,2$,
and the Gaussian initial state $x_0 = \exp(-q_1^2-q_2^2)$.
This example is chosen for the following reasons.
First, the exact solution repeats with the period of $T_p=20$, hence we can estimate the error of our scheme as the difference between the solution after a number of periods and the initial state, $\|x(T)-x_0\|$ for $T=20m$, $m \in \mathbb{N}$.
Second, \eqref{eq:conv} possesses both types of invariants: the solution mass $c^*x = c^*x_0$, where $c=(1,\ldots,1)^\top$, and the second norm, $\|x\|_2 = \|x_0\|_2$.
Third, the discrete solution of the pure convection is prone to developing spurious oscillations when the discretization is not accurate enough.
For the central difference scheme, this requires taking rather fine grids, with $n$ ranging from $1024$ to $4096$,
which makes the problem large enough to apply the tensor decompositions.

In order to increase the efficiency of the TT methods, we apply them to the \emph{quantized} tensors \cite{khor-qtt-2011}:
instead of separating just two indices in $X(i_1,i_2)$,
we split each index into its binary digits,
$$
i_k = i_{k,1} + 2 (i_{k,2}-1) + \cdots + 2^{L-1} (i_{k,L}-1), \quad i_{k,l} \in \{1,2\},
$$
and consider the solution as a $2L$-dimensional tensor, $\x(i_{1,1},\ldots,i_{2,L})$.
Now the TT decomposition can reduce the storage cost down to $\mathcal{O}(r^2L) = \mathcal{O}(r^2\log n )$, in contrast to $\mathcal{O}(rn)$ in the low-rank decomposition of $X(i_1,i_2)$ or $n^2$ in the full representation of $x$.
The ODE matrix can be constructed in this quantized TT representation exactly \cite{khkaz-conv-2013}.

%

%

\hypersetup{pdfborder={0 0 0}}
\begin{figure}[t]
\centering
\caption{Convection example. Degeneracy of $\|u\|_2$ (solid lines) and $c^*u$ (dashed lines) vs. time for Chebyshev (left) and Crank-Nicolson (right) schemes.}
\label{fig:conv_dobs}
\resizebox{0.49\linewidth}{!}{\begin{tikzpicture}
\begin{axis}[%
  xmode=normal,ymode=log,
  cycle list name=bytwo,
  xmin=0, xmax=5,
  xtick={1,2,3,4,5},
  xticklabels={1,2,3,4,~},
  ymax=2e-9,
  y label style={at={(-0.03,1.0)},anchor=south west},
  xlabel={$t/T_p$},
  ylabel={error},
  yminorticks=true,
  legend style={at={(0.01,0.99)},anchor=north west}]

\addplot+[] table[header=true, x=t, y=e24]{dats/conv_tamen_dobs_new.dat}; \label{dobs_48};
\addplot+[] table[header=true, x=t, y=e14]{dats/conv_tamen_dobs_new.dat};
\addplot+[] table[header=true, x=t, y=e25]{dats/conv_tamen_dobs_new.dat}; \label{dobs_58};
\addplot+[] table[header=true, x=t, y=e15]{dats/conv_tamen_dobs_new.dat};
\addplot+[] table[header=true, x=t, y=e25_16]{dats/conv_tamen_dobs_new.dat}; \label{dobs_516};
\addplot+[] table[header=true, x=t, y=e15_16]{dats/conv_tamen_dobs_new.dat};

\node [fill=none,anchor=south west] (ndobs_516) at (axis cs: 0.0,3e-11) {~\ref{dobs_516}~$\eps=10^{-5},~\I=16$};
\node [fill=none,anchor=south,above=0cm of ndobs_516] (ndobs_58)  {\ref{dobs_58}~$\eps=10^{-5},~\I=8$};
\node [fill=none,anchor=south,above=0cm of ndobs_58]  {\ref{dobs_48}~$\eps=10^{-4},~\I=8$};

\end{axis}
\end{tikzpicture}%
}
\resizebox{0.49\linewidth}{!}{\begin{tikzpicture}
\begin{axis}[%
  xmode=normal,ymode=log,
  cycle list name=bytwo,
  xmin=0, xmax=5,
  xtick={1,2,3,4,5},
  xticklabels={1,2,3,4,~},
  ymax=2e-9,
  y label style={at={(-0.03,1.0)},anchor=south west},
  xlabel={$t/T_p$},
  ylabel={error},
  yminorticks=true,
  legend style={at={(0.01,0.99)},anchor=north west}]

\addplot+[] table[header=true, x=t, y=e2_4_1025]{dats/conv_cn_dobs.dat}; \label{dobs_4};
\addplot+[] table[header=true, x=t, y=e1_4_1025]{dats/conv_cn_dobs.dat};
\addplot+[] table[header=true, x=t, y=e2_1025]{dats/conv_cn_dobs.dat}; \label{dobs_5};
\addplot+[] table[header=true, x=t, y=e1_1025]{dats/conv_cn_dobs.dat};
\addplot+[] table[header=true, x=t, y=e2_4097]{dats/conv_cn_dobs.dat}; \label{dobs_4097};
\addplot+[] table[header=true, x=t, y=e1_4097]{dats/conv_cn_dobs.dat};

\node [fill=none,anchor=south west] (ndobs_4097) at (axis cs: 0.0,3e-11) {~\ref{dobs_4097}~$\eps=10^{-5},~\I=4097$};
\node [fill=none,anchor=south,above=0cm of ndobs_4097] (ndobs_5)  {\ref{dobs_5}~$\eps=10^{-5},~\I=1025$};
\node [fill=none,anchor=south,above=0cm of ndobs_5]  {\ref{dobs_4}~$\eps=10^{-4},~\I=1025$};

\end{axis}
\end{tikzpicture}%
}
\end{figure}
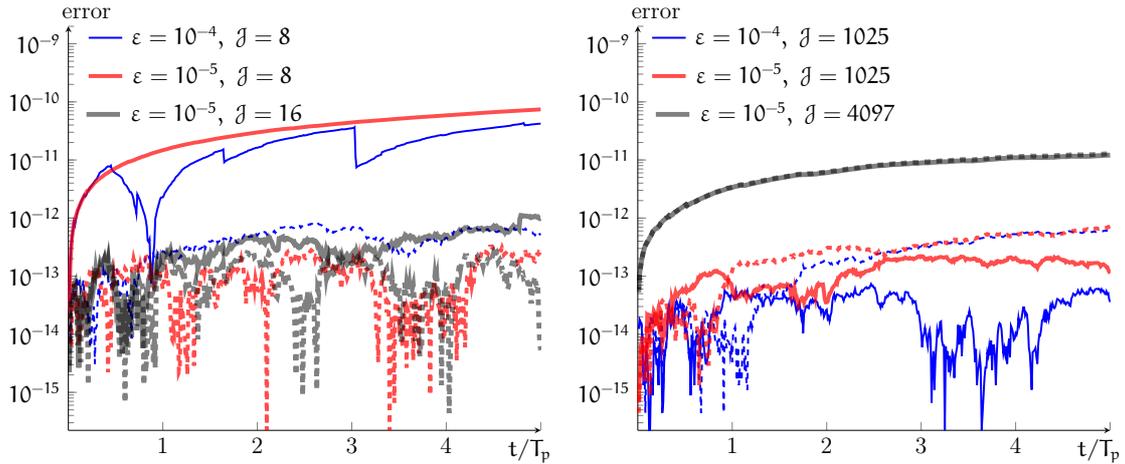
\hypersetup{pdfborder={0 0 0.1}}

First, we confirm conservation of the invariants.
We fix the time interval splitting to $h_{\ell}=0.2$ for all $\ell=1,\ldots,T/h_{\ell}$, the spatial grid size $n=4096$, and vary the number of Chebyshev or Crank-Nicolson points, as well as the accuracy threshold $\eps$.
In Fig. \ref{fig:conv_dobs} we show how the errors in both conserving quantities evolve with time.
We observe that the error in the invariants is much smaller than the tensor approximation threshold in all cases.
However, insufficient number of Chebyshev points can increase the error in the second norm (Fig. \ref{fig:conv_dobs}, left).
The Crank-Nicolson scheme preserves both invariants up to the machine precision for any number of points.
In fact, it manifests the opposite situation that the errors are larger for $\I=4097$ points due to a larger condition number of the matrix in \eqref{eq:xt-system}.
This shows that, although the direct cost of the tensor schemes depends mildly on the grid sizes,
it is still recommended to avoid too fine grids due to the conditioning issues.

\begin{figure}[t]
\centering
\caption{Convection example, TT ranks for Chebyshev (left) and Crank-Nicolson (right) schemes.}
\label{fig:conv_rnk}
\resizebox{0.49\linewidth}{!}{\begin{tikzpicture}
\begin{axis}[%
  xmode=normal,ymode=normal,
  cycle list name=eta,
  xtick={1,2,3,4,5},
  xticklabels={1,2,3,4,~},
  ymin=15, ymax=38,
  xlabel={$t/T_p$},
  ylabel={max TT rank},
  yminorticks=true,
  legend style={at={(0.01,0.99)},anchor=north west}]

  \addplot+[no marks] table[header=true, x=t, y=r_4_8] {dats/conv_tamen_rnk.dat}; \addlegendentry{$\eps=10^{-4}$, $\I=8$};
  \addplot+[no marks] table[header=true, x=t, y=r_5_8] {dats/conv_tamen_rnk.dat}; \addlegendentry{$\eps=10^{-5}$, $\I=8$};
  \addplot+[no marks] table[header=true, x=t, y=r_5_16] {dats/conv_tamen_rnk.dat}; \addlegendentry{$\eps=10^{-5}$, $\I=16$};

%
%
%
%
%
%
%

\end{axis}
\end{tikzpicture}%
}
\resizebox{0.49\linewidth}{!}{\begin{tikzpicture}
\begin{axis}[%
  xmode=normal,ymode=normal,
  cycle list name=eta,
  xtick={1,2,3,4,5},
  xticklabels={1,2,3,4,~},
  ymin=15, ymax=38,
  xlabel={$t/T_p$},
  ylabel={max TT rank},
  yminorticks=true,
  legend style={at={(0.01,0.99)},anchor=north west}]

  \addplot+[no marks] table[header=true, x=t, y=r_4_1025] {dats/conv_cn_rnk.dat}; \addlegendentry{$\eps=10^{-4}$, $\I=1025$};
  \addplot+[no marks] table[header=true, x=t, y=r_5_1025] {dats/conv_cn_rnk.dat}; \addlegendentry{$\eps=10^{-5}$, $\I=1025$};
  \addplot+[no marks] table[header=true, x=t, y=r_5_4097] {dats/conv_cn_rnk.dat}; \addlegendentry{$\eps=10^{-5}$, $\I=4097$};

%
%
%
%
%
%
%

\end{axis}
\end{tikzpicture}%
}
\end{figure}
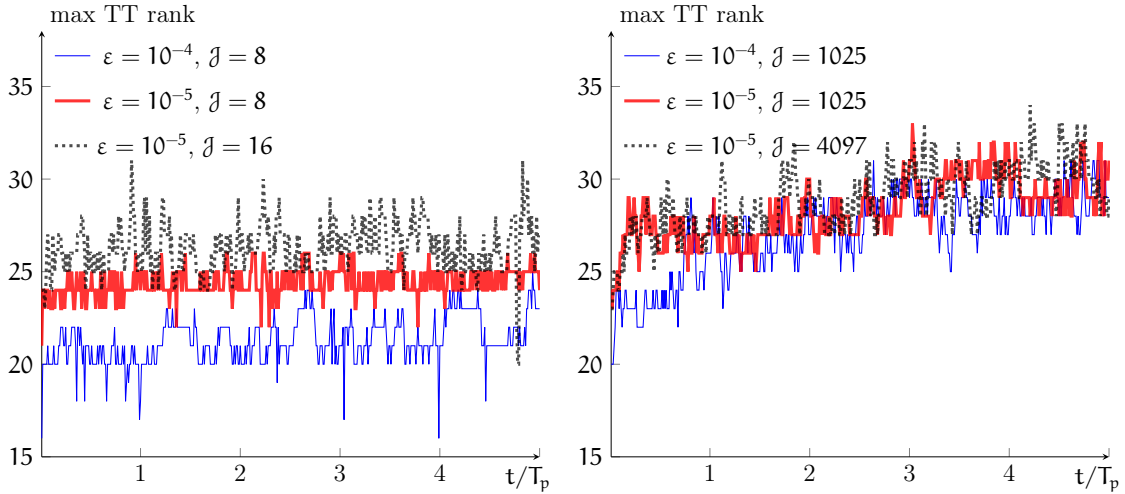

\begin{figure}[t]
\centering
\caption{Convection example, CPU times in seconds (left) and time intervals found in the adaptive regime (right) with Chebyshev and Crank-Nicolson (CN) schemes.}
\label{fig:conv_ttimes}
\resizebox{0.49\linewidth}{!}{\begin{tikzpicture}
\begin{axis}[%
  xmode=normal,ymode=normal,
  cycle list name=eta,
  xtick={1,2,3,4,5},
  xticklabels={1,2,3,4,~},
  xlabel={$t/T_p$},
  ylabel={cumulative CPU time},
  yminorticks=true,
  legend style={at={(0.99,0.01)},anchor=south east}]


\addplot+[] table[header=true, x=t, y=T_cheb_8]{dats/conv_tamen-cn_ttimes_cum.dat}; \addlegendentry{Cheb, $\I=8$};
\addplot+[] table[header=true, x=t, y=T_cheb_16]{dats/conv_tamen-cn_ttimes_cum.dat}; \addlegendentry{Cheb, $\I=16$};
\addplot+[] table[header=true, x=t, y=T_cn_513]{dats/conv_tamen-cn_ttimes_cum.dat}; \addlegendentry{CN, $\I=513$};
\addplot+[] table[header=true, x=t, y=T_cn_2049]{dats/conv_tamen-cn_ttimes_cum.dat}; \addlegendentry{CN, $\I=2049$};


\end{axis}
\end{tikzpicture}%
}
\resizebox{0.49\linewidth}{!}{\begin{tikzpicture}%
\begin{axis}[%
     xmode=normal,ymode=normal,
     cycle list name=eta,
     xtick={1,2,3,4,5},
     xticklabels={1,2,3,4,~},
     ymin=0,ymax=0.27,
     ytick={0.05,0.1,0.15,0.2,0.25},
     yticklabels={0.05,0.1,0.15,0.2,0.25},
     xlabel={$t/T_p$},
     ylabel={$h_{\ell}$},
     legend style={at={(0.99,0.01)},anchor=south east}]
  \addplot+[] table[header=false,x index=0,y index=1] {dats/conv_tamen_adapt_dt.dat}; \addlegendentry{Cheb, $\I=8$};
  \pgfplotsset{cycle list shift=1}
  \addplot+[] table[header=false,x index=0,y index=1] {dats/conv_cn_adapt_dt.dat}; \addlegendentry{CN, $\I=513$};
\end{axis}
\end{tikzpicture}%
}
\end{figure}
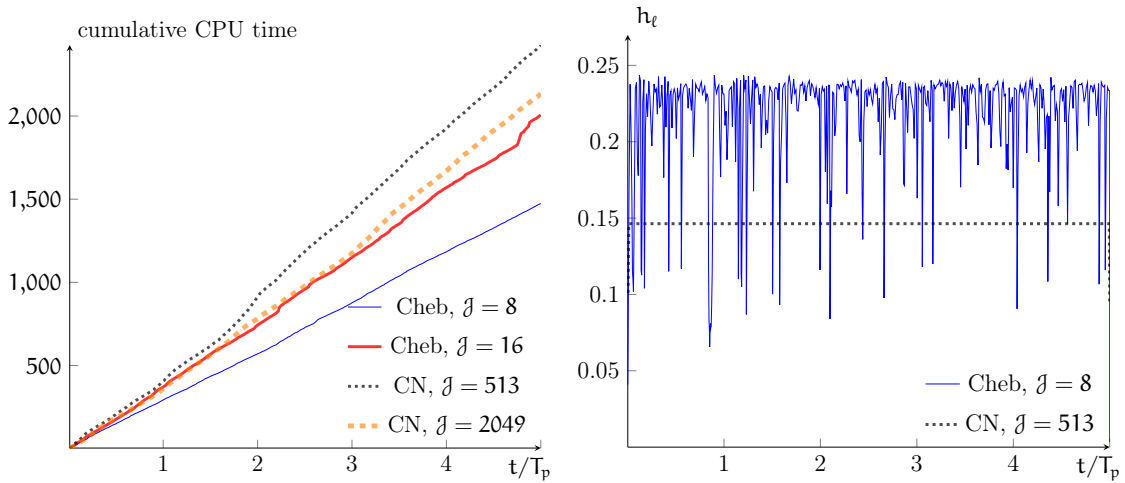

The evolution of TT ranks with time is shown in Fig. \ref{fig:conv_rnk}.
Since the shape of the exact solution remains unchanged, its ranks should be the same for all times.
We see that the ranks are indeed stable with time, particularly with the Chebyshev scheme.
In the Crank-Nicolson scheme, the ranks grow slightly towards the end of the 5-period interval.
This is also reflected by a slightly larger CPU time, see Fig. \ref{fig:conv_ttimes} (left).

\begin{table}[t]
\centering
\caption{Convection, CPU times (seconds) and errors for different time interval lengths.}
\label{tab:conv_ttimes_tau}
\begin{tabular}{c|ccc|c|ccc|c}
Scheme                      & \multicolumn{4}{c|}{Chebyshev, $\I=8$} &   \multicolumn{4}{c}{Crank-Nicolson, $\I=513$}  \\
Time step                   & $0.1$   & $0.2$  & $0.4$  & $100$      &  $0.1$ & $0.2$  & $0.4$  & $100$  \\ \hline
CPU time                    & 1098.8  & 1391.1 & 5069.9 & 2014.7     & 4326.0 & 2213.7 & 10246.1 & 2652.9 \\
$\frac{10^4\cdot \|x-x_\star\|}{\|x_\star\|}$ & 0.85 & 6.22 & 5.65 & 2.81  & 3.33   & 3.52   & 2.22   & 2.62 \\
\end{tabular}
\end{table}

Now we consider how the tAMEn algorithm depends on the time interval splitting.
In Table \ref{tab:conv_ttimes_tau} we show the CPU times and the errors of $x(T)$ with respect to the reference solution $x_\star$, computed with the Chebyshev scheme with $\I=16$ points on $h_{\ell}=0.2$ and $\eps=10^{-7}$.
For small time steps ($h_{\ell}=0.1,0.2$ and $0.4$), we turn the adaptation off.
However, we also start from the entire interval $100$ and let the algorithm split it automatically.
Due to rejections of some time steps, the CPU time of the adaptive method is larger than the cost of the optimal splitting ($h_{\ell}=0.1$ for Chebyshev and $h_{\ell}=0.2$ for Crank-Nicolson schemes),
but the overhead never exceeds a factor of $2$.
Moreover, the adaptive algorithm is faster than the non-adaptive one with improperly chosen time steps.
Fig. \ref{fig:conv_ttimes} (right) shows the time steps determined by the adaptive method.
We see that the average step lies between $0.1$ and $0.2$.
Interestingly, the low-order Crank-Nicolson scheme is more robust in estimating the error, and hence the time step.

\begin{table}[t]
\centering
\caption{Convection example. CPU times (seconds) and errors in different methods and parameters for time splitting $h_{\ell}=0.2$.}
\label{tab:conv_timeerror}
\begin{tabular}{c|cc|cc|cc}
                               & \multicolumn{2}{c|}{tAMEn}  & \multicolumn{2}{c|}{KSL} & \multicolumn{2}{c}{Full CN} \\
                               & Cheb, $\I=8$ & CN, $\I=513$ & $\I=16$ & $\I=512$ & $\I=16$ & $\I=64$    \\ \hline
CPU time                       & 1391.1          & 2213.7    & 694.9   & 14159    & 170732  & 102294   \\
$\frac{10^3 \cdot \|x(T)-x_0\|}{\|x_0\|}$ & 2.36 & 2.16      & 583.8   & 7.44     & 16.0    & 2.93     \\
\end{tabular}
\end{table}

%
%
%

Finally, we benchmark tAMEn against the standard Crank-Nicolson method without the TT decomposition
and the Riemannian TT time integrator \cite{lob-ttksl-2015}.
We split the time into intervals of length $h_{\ell}=0.2$,
but each interval is further partitioned into $\I$ individual time steps,
on which the full Crank-Nicolson or Riemannian integration is carried out.
The Riemannian integrator projects the dynamical equations directly onto the Riemannian manifold of the TT representation, using the so-called \emph{Dirac-Frenkel} principle \cite{lubich-koch-dynten-2010}.
The projected equations can be split with respect to the different TT blocks and solved subsequently, using the so-called KSL propagator \cite{lob-ttksl-2015}\footnote{The multi-dimensional Matlab version \texttt{tt\_ksl\_ml.m} was implemented by the author in collaboration with I. Oseledets, and is available within \href{http://github.com/oseledets/TT-Toolbox}{TT-Toolbox}.}
This scheme works with the TT decomposition of only one snapshot at a time, which requires smaller TT ranks.
However, it requires integrating backward in time, which can introduce numerical instabilities for large time steps.
In Table \ref{tab:conv_timeerror} we see that for $\I=16$ the solution becomes qualitatively incorrect.
For a smaller time step the scheme is stable, but a large number of time steps leads to a large computational time.
The full Crank-Nicolson method is even slower, since each time step is more expensive.
In fact, the CPU time is larger for smaller number of time steps.
This is due to a larger condition number of the matrix in the implicit step.

\subsection{Chemical master equation}
In the second experiment, we investigate an example with a steady state,
the chemical master equation (CME), describing stochastic kinetics model of the $\lambda$-phage virus \cite{hegland-cme-2007,jahnke-cme-2008,dkh-cme-2014}.
Using the Finite State Projection \cite{munsky-fsp-2006}, the CME is turned into a large-scale ODE,
\begin{equation}
\dfrac{d x}{d t} = A x, \qquad A = \sum_{m=1}^M \left(J^{z_1^m} \otimes \cdots \otimes J^{z_d^m} - I\right) \mathrm{diag}(w^m),
\label{eq:cme}
\end{equation}
Here, $J^z$ is the order-$z$ shift matrix, defined as follows: $J^0=I$, $J^1 = \mathrm{tridiag}(1,0,0)$, $J^z = (J^1)^z$ for $z > 1$, and $J^z=(J^{-z})^\top$ for $z<0$.
The vector $\mathbf{z}^m = (z^m_1,\ldots,z^m_d)$ is the so-called \emph{stoichiometric} vector,
$w^m = w^m(i_1,\ldots,i_d)$ is the \emph{propensity} rate of the $m$-th reaction, and $\mathrm{diag}(w^m)$ constructs a $N \times N$ diagonal matrix from all elements of $w^m$.
The total size of the problem is $N=\prod_{k=1}^d n_k$, since each index is assumed to vary in the range $i_k = 0,\ldots,n_k-1$.
The indices $i_1,\ldots,i_d$ denote the so-called \emph{copy numbers} (numbers of molecules) of $d$ reacting species (e.g. proteins), and
the solution $\x(i_1,\ldots,i_d,t)$ is the distribution function, which defines the probability that at the time $t$, the system contains $i_1$ molecules of the first protein, $i_2$ of the second, and so on.

The particular $\lambda$-phage model contains $d=5$ species and $M=10$ reactions.
The stoichiometric vectors and propensities are given in Table \ref{tab:lphage} ($\mathbf{e}_1,\ldots,\mathbf{e}_5$ are unit vectors of size $5$).

\begin{table}[t]
\centering
\caption{Reactions in the $\lambda$-phage model.}
\label{tab:lphage}
\begin{tabular}{l|ll|ll}
       & \multicolumn{2}{c|}{Generation} & \multicolumn{2}{c}{Destruction} \\[0.5em]\hline
$S_1$  & $w^{1} = \dfrac{0.06}{0.12+i_2}$, & $\mathbf{z}^{1} = \mathbf{e}_1$ & $w^{2} = 0.0025 \cdot i_1$, & $\mathbf{z}^{2} = -\mathbf{e}_1$ \\[0.5em]
$S_2$  & $w^3 = \dfrac{(1+i_5)\cdot 0.6}{0.6 + i_1}$, & $\mathbf{z}^{3} = \mathbf{e}_2$ & $w^4 = 0.0007 \cdot i_2$, & $\mathbf{z}^{4}= -\mathbf{e}_2$ \\[0.5em]
$S_3$  & $w^5=\dfrac{0.15 \cdot i_2}{i_2+1}$, & $\mathbf{z}^{5}= \mathbf{e}_3$ & $w^6=0.0231 \cdot i_3$, & $\mathbf{z}^{6}= -\mathbf{e}_3$ \\[0.5em]
$S_4$  & $w^7=\dfrac{0.3 \cdot i_3}{i_3+1}$, & $\mathbf{z}^{7}= \mathbf{e}_4$ & $w^8=0.01 \cdot i_4$, & $\mathbf{z}^{8}= -\mathbf{e}_4$ \\[0.5em]
$S_5$  & $w^9=\dfrac{0.3 \cdot i_3}{i_3+1}$, & $\mathbf{z}^{9}= \mathbf{e}_5$ & $w^{10}=0.01 \cdot i_5$, & $\mathbf{z}^{10}= -\mathbf{e}_5$ \\[0.5em]
\end{tabular}
\end{table}

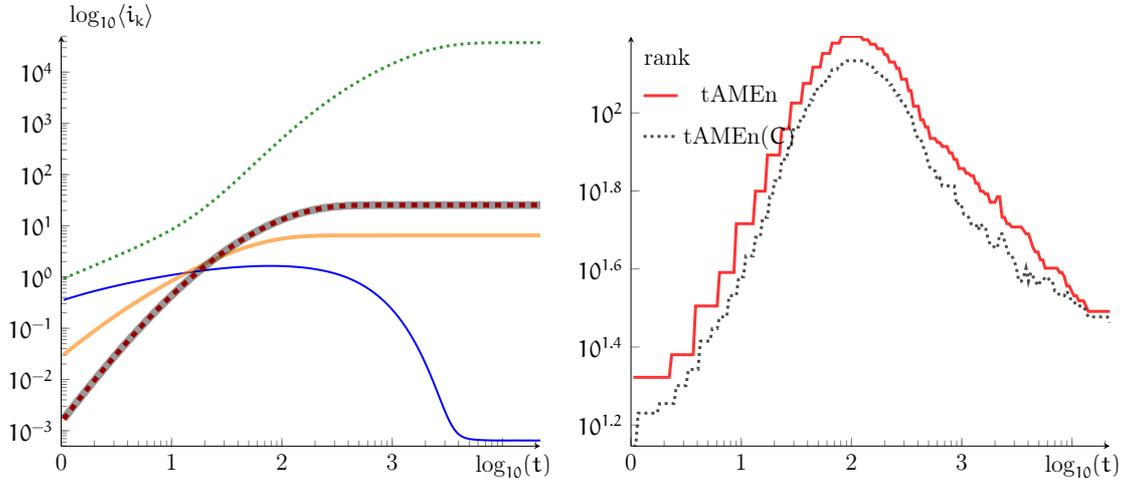
\begin{figure}[t]
\centering
\caption{CME example, $\langle i_k \rangle$ (left) and maximal TT ranks in tAMEn with and without $C$-enrichment (right)}
\label{fig:cme-meanconc}
\resizebox{0.49\linewidth}{!}{%
\begin{tikzpicture}%
\begin{axis}[%
xmode=log,ymode=log,
xmin=1, xmax=2.2e4,
xtick={1e0,1e1,1e2,1e3,1e4},
xticklabels={0,1,2,3,~},
xlabel={$\log_{10}(t)$},
ylabel={$\log_{10}\langle i_k \rangle$},
ymin=5e-4, ymax=5e4,
legend style={at={(1.0,0.055)},anchor=south east}]

\addplot+[no marks,line width=1.0pt,blue,opacity=1.0] table[header=false,x index=0,y index=1]{dats/lphage-meanconc-new.dat};
\addplot+[no marks,line width=1.5pt,black!50!green,opacity=0.8,dotted] table[header=false,x index=0,y index=2]{dats/lphage-meanconc-new.dat};
\addplot+[no marks,line width=2.0pt,orange,opacity=0.6] table[header=false,x index=0,y index=3]{dats/lphage-meanconc-new.dat};
\addplot+[no marks,line width=2.5pt,red,opacity=1.0,dashed] table[header=false,x index=0,y index=4]{dats/lphage-meanconc-new.dat};
\addplot+[no marks,line width=4.0pt,black,opacity=0.4,solid] table[header=false,x index=0,y index=5]{dats/lphage-meanconc-new.dat};
\end{axis}
\end{tikzpicture}%
}
\resizebox{0.49\linewidth}{!}{\begin{tikzpicture}
\begin{axis}[%
  cycle list name=eta,
  xmode=log,ymode=log,
  xmin=1, xmax=2.2e4,
  xtick={1e0,1e1,1e2,1e3,1e4},
  xticklabels={0,1,2,3,~},
  xlabel={$\log_{10}(t)$},
  yminorticks=true,
  legend style={at={(0.01,0.9)},anchor=north west}]

  \pgfplotsset{cycle list shift=1};
  \addplot+[] table[header=false,x index=0,y index=1]{dats/lphage-rnk.dat}; \addlegendentry{tAMEn};
  \addplot+[] table[header=false,x index=0,y index=2]{dats/lphage-rnk.dat}; \addlegendentry{tAMEn($C$)};


   
\node[anchor=north west] at (rel axis cs: 0.01,0.99) {rank};

\end{axis}
\end{tikzpicture}%
}
\end{figure}
\hypersetup{pdfborder={0 0 0.1}}

As the initial state, we choose all-zero copy numbers with probability $1$, i.e. $x_0(i_1,\ldots,i_5)=1$ when $i_1=\cdots=i_5=0$, and $0$, otherwise.
Under certain conditions \cite{khammash-ergodicity-2014}, fulfilled for the $\lambda$-phage model, and infinite ranges of $i_k$, the CME \eqref{eq:cme} converges to a unique stationary state $x_{\infty}$.
For practical computations, we truncate the state space to $n_1 \times \cdots \times n_5 = 128 \times 65536 \times 64 \times 64 \times 64$, respectively, since the probability outside this box is negligible.
In order to preserve existence of the stationary state \cite{jahnke-cme-red-2011},
we adjust the propensities of the generation reactions such that
$$
w^{2k-1}(i_1,\ldots,i_d) = 0 \quad \mbox{if} \quad i_k = n_k-1, \quad k=1,\ldots,d.
$$
This also guarantees that $A^*e=0$, therefore the probability \emph{normalization} $e^* x=1$ is conserved.

The statistical outputs of interest are the \emph{mean copy numbers},
\begin{equation}
\langle i_k \rangle = \frac{\mathtt{i}_k^* \psi}{e^* \psi}, \quad \mathtt{i}_k = \mathtt{e}^{(1)} \otimes \cdots \otimes \mathtt{e}^{(k-1)} \otimes \{i_k\} \otimes \mathtt{e}^{(k+1)} \otimes \cdots \otimes \mathtt{e}^{(d)} \in \mathbb{R}^N,
\label{eq:mean_i}
\end{equation}
where $\mathtt{e}^{(p)}$ are the all-ones vectors of size $n_p$.
In order to preserve the normalization, we add the vector of ones to the enrichment \eqref{eq:enrich_both}.
However, we can also keep the quantities of interest in the TT representation in order to make statistics more accurate.
Therefore, we use $6$ enrichment columns, $C = \begin{bmatrix}e & \mathtt{i}_1 & \cdots & \mathtt{i}_5\end{bmatrix}$.
The Quantized TT representation of $C$ has TT ranks up to $6$, and the ranks of the residual \eqref{eq:zk} are set to $1$.
We compare tAMEn implementations with and without the additional enrichment by $C$.
For the fair comparison, we set the residual ranks equal to those of $C$ plus $1$ in the version without the $C$-enrichment.

The tAMEn algorithm is run in the fully adaptive regime, applied to the time interval $T=22000$.
We estimate the errors directly in the quantities of interest
by taking the log average values of the mean copy numbers at certain times,
\begin{equation}
\mathcal{E}_{\langle i \rangle}(t) = \exp\left(\frac{1}{5}\sum_{k=1}^{5} \log \frac{\left|\langle i_k(t) \rangle - \langle i_k^\star(t) \rangle\right|}{\langle i_k^\star(t) \rangle} \right)
\label{eq:cmeerr}
\end{equation}
We vary the accuracy thresholds $\eps$ from $10^{-2}$ to $10^{-5}$, and use the values computed with $\eps=3 \cdot 10^{-7}$ as the reference $\langle i_k^\star \rangle$.
In addition to the two versions of tAMEn, we present the results of the classical Stochastic Simulation Algorithm (SSA) \cite{gillespie-ssa-1976} for comparison.

\begin{figure}[htb]
\caption{CME example, errors \eqref{eq:cmeerr} in the mean copy numbers for $t=2000$ (left) and $t=22000$ (right) versus computational Work (CPU time) for SSA and tAMEn with and without the $C$-enrichment.}
\label{fig:cmee}
\resizebox{0.49\linewidth}{!}{
\begin{tikzpicture}
\begin{axis}[%
  xmode=log,ymode=log,
  cycle list name=eta,
  y label style={at={(-0.03,1.0)},anchor=south west},
  xlabel={CPU time},
  ylabel={$\mathcal{E}_{\langle i \rangle}$},
  xmax={3e4},
  yminorticks=true,
  legend style={at={(0.01,0.01)},anchor=south west}]

\addplot+[] coordinates{(5.8278e+02, 1.3284e-02)
                        (1.0900e+03, 1.1069e-02)
                        (2.1426e+03, 5.0672e-03)
                        (4.1099e+03, 6.8354e-03)
                        (8.8998e+03, 2.4563e-03)
                        (1.8370e+04, 2.4045e-03)
                        }; \addlegendentry{SSA};

\addplot+[] coordinates{(2.0252e+02, 4.7816e-03)
                        (4.3729e+02, 1.2079e-03)
                        (8.0832e+02, 5.7301e-04)
                        (1.6524e+03, 4.2221e-04)
                        (3.1832e+03, 1.2268e-04)
                        (6.6883e+03, 3.7919e-05)
                        (1.3085e+04, 1.0194e-05)
                        }; \addlegendentry{tAMEn};

\addplot+[] coordinates{(3.5559e+02, 8.0999e-04)
                        (6.0765e+02, 4.1398e-04)
                        (9.4883e+02, 1.8097e-04)
                        (1.9745e+03, 2.6332e-05)
                        (3.7236e+03, 1.3350e-05)
                        (7.3754e+03, 3.2558e-06)
                        (1.3449e+04, 2.1464e-06)
                        }; \addlegendentry{tAMEn($C$)};

\addplot+[domain=5e2:2e4,no marks] {sqrt(1/x)*0.3}; \addlegendentry{$W^{-0.5}$};
\addplot+[domain=5e2:2e4,no marks,purple,dashed] {1e2*(x^(-1.7))}; \addlegendentry{$W^{-1.7}$};
\end{axis}
\end{tikzpicture}%
}
\resizebox{0.49\linewidth}{!}{
\begin{tikzpicture}
\begin{axis}[%
  xmode=log,ymode=log,
  cycle list name=eta,
  y label style={at={(-0.03,1.0)},anchor=south west},
  xlabel={CPU time},
  ylabel={$\mathcal{E}_{\langle i \rangle}$},
  xmax={3e4},
  yminorticks=true,
  legend style={at={(0.01,0.01)},anchor=south west}]

\addplot+[] coordinates{(5.8278e+02, 1.1302e-02)
                        (1.0900e+03, 1.1678e-02)
                        (2.1426e+03, 6.5641e-03)
                        (4.1099e+03, 5.1346e-03)
                        (8.8998e+03, 3.5369e-03)
                        (1.8370e+04, 3.0537e-03)
                        }; \addlegendentry{SSA};

\addplot+[] coordinates{(2.0252e+02, 1.0603e-03)
                        (4.3729e+02, 5.1729e-04)
                        (8.0832e+02, 4.1026e-04)
                        (1.6524e+03, 1.8410e-04)
                        (3.1832e+03, 1.9013e-05)
                        (6.6883e+03, 1.5313e-05)
                        (1.3085e+04, 7.6647e-06)
                        }; \addlegendentry{tAMEn};

\addplot+[] coordinates{(3.5559e+02, 1.5811e-03)
                        (6.0765e+02, 6.4921e-04)
                        (9.4883e+02, 3.5828e-04)
                        (1.9745e+03, 7.5014e-05)
                        (3.7236e+03, 2.0972e-05)
                        (7.3754e+03, 8.1294e-06)
                        (1.3449e+04, 2.1227e-06)
                        }; \addlegendentry{tAMEn($C$)};

\addplot+[domain=5e2:2e4,no marks] {sqrt(1/x)*0.3}; \addlegendentry{$W^{-0.5}$};
\addplot+[domain=5e2:2e4,no marks,purple,dashed] {33*(x^(-1.7))}; \addlegendentry{$W^{-1.7}$};
\end{axis}
\end{tikzpicture}%
}
\end{figure}
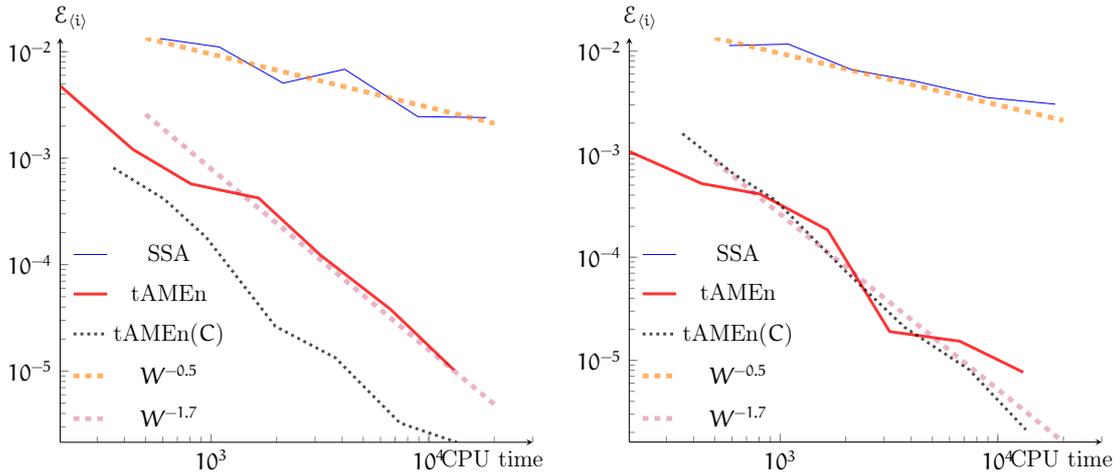

In Fig. \ref{fig:cme-meanconc} we show the evolution of the mean copy numbers and TT ranks.
Interestingly, the ranks with the $C$-enrichment are even smaller, since the specialized frame matrices constitute better bases for the solution.
The computational times and errors are shown in Fig. \ref{fig:cmee}.
We see that the normalization-preserving solution is systematically more efficient in terms of the cost/accuracy ratio, compared to the residual-only enrichment.
Moreover, the direct solution of the CME in the TT format is much faster than the stochastic simulation,
since large times and copy numbers require a large number of trajectories and time steps in SSA.

\section{Conclusion}
We have proposed and studied an alternating iterative algorithm for approximate solution of ordinary differential equations in the TT format.
The method combines advances of DMRG techniques and classical iterative methods of linear algebra.
Started from the solution at the previous time interval as the initial guess, it often converges in 2---4 iterations,
and delivers accurate solution even for strongly non-symmetric matrices in the right-hand side of an ODE.
%
The numerical experiments reveal a promising potential of this method in long time simulations when the solution admits a low-rank decomposition.
For example, nuclear magnetic resonance models
can be approached directly, without any a priori reduction of the original Hilbert space \cite{sdwk-nmr-2014}.
On the other hand, it was important for tAMEn algorithm to split the time integration into steps.
It remains an open question how this method would perform for nonlinear or inverse problems, where the solution must be computed on the solid time interval.

\end{document}